\documentclass[12pt]{iopart}
\usepackage{iopams}
\usepackage{algorithm}
\usepackage{algpseudocode}
\usepackage{subfigure}

\expandafter\let\csname equation*\endcsname\relax

\expandafter\let\csname endequation*\endcsname\relax

\usepackage{graphicx}
\usepackage{amsmath}
\usepackage{amsthm}
\usepackage{amssymb}
\usepackage{graphicx}
\usepackage{color}
\usepackage{tikz}
\usepackage{enumitem}

\begin{document}
\title[Well-posed Questions for Ill-posed Problems]{Well-posed Questions for Ill-posed Inverse Problems: a Note in Memory of Pierre Sabatier}
\author{Gaoming Chen$^1$, Fadil Santosa$^2$, William W. Symes$^3$}
\address{$^1$ $^2$Department of Applied Mathematics and Statistics, Johns Hopkins University, Baltimore, MD 21218, USA}
\address{$^3$Department of Computational Applied Mathematics and Operations Research, Rice University, Houston, TX 77251}
\ead{fsantos9@jhu.edu}

\vspace{10pt}
\begin{indented}
\item[]March 2025
\end{indented}

\maketitle

\begin{abstract}
    Professor Pierre Sabatier contributed much to the study of inverse problems in theory and practice. Two of these contributions were a focus on theory that actually supports practice, and the identification of well-posed aspects of inverse problems that may quite ill-posed. This paper illustrates these two themes in the context of Electrical Impedance Tomography (EIT), which is both very ill-posed and very practical. We show that for a highly constrained version of this inverse problem, in which a small elliptical inclusion in a homogeneous background is to be identified, optimization of the experimental design (that is, electrode locations) vastly improves the stability of the solution.
\end{abstract}

\section{Introduction}

Pierre Sabatier founded the journal Inverse Problems, which 40 years on continues to serve its community and the field of inverse problems with distinction. In a very real sense, he also founded the field, by bringing together conceptual threads and nomenclature that were current in various parts of science and mathematics but not yet in fullest possible contact. His ideas and terminology subsequently pervaded the subject and continue to do so.

This brief communication centers around two ideas that Sabatier emphasized many times in his copious writings on inverse problems: the quest for
``well-posed questions for ill-posed problems'' \cite{SabatierSIAM:84}, and the intrinsically interdisciplinary nature of the subject \cite{SabatierJMP:00}. We will illustrate these in the context of Electrical Impedance Tomography (EIT), an inverse problem in electrostatics with many biomedical and engineering applications \cite{adler, borcea}. EIT seeks to infer the conductivity distribution in the interior of a region from the measured voltage drops between pairs of electrodes on the boundary in response to specified current flow between them. This problem is spectacularly ill-posed. For a strongly constrained version, we ask how the electrodes should be distributed to minimize a measure of solution uncertainty. This experimental design question turns out to have a computable and relatively stable answer. It is a well-posed question for an ill-posed problem.   

The concept of well-posed questions is also found in the writings of D. D. Jackson \cite{ddjackson:72}, for example. The phrase encapsulates a recognition that inverse problems in science and engineering tend to be ill-posed, i.e., the mapping between model and data (a framework that Sabatier did much to popularize) may map very different models to nearby data. Thus the model (solution of the inverse problem) is poorly determined, if one regards the inverse problem as a functional equation (or more generally as a fit optimization problem). In particular, small uncertainty in the data leads to large uncertainty in the model, in such formulations. Even for such ill-posed problems, however, it may be possible to identify functions of the solution, the values of which have small uncertainty with small data uncertainty. The question, whether such function values lie in specified intervals given assumed data error bounds, then has a definite answer, yes or no. That is, such questions are examples of well-posed questions about ill-posed problems.

Sabatier frequently contended that the study of inverse problems is inherently interdisciplinary, or applied. He rejected as uninteresting purely mathematical exercises without any implications for science and engineering practice \cite{SabatierJMP:00}. His own interpretation of this dictum was creatively subtle. He and his collaborators contributed many indisputably practical results, in geophysics and other areas. However he also maintained interests in many areas, such as 1D models of quantum inverse scattering, that could be termed proto-practical.

EIT is very much an applied inverse problem, with laboratory implementations in biomedical imaging and nondestructive testing \cite{adler, borcea}. It has also received a great deal of mathematical attention, beginning with the seminal work of Calderon \cite{calderon}. As the EIT problem is strongly ill-posed, we restrict our discussion to a heavily constrained 2D special case. The conductivity model to be determined consists of a small elliptical inclusion in a disk of known conductivity. Such an ellipse depends on five parameters (center coordinates, axes, and orientation angle). The conductivity anomaly in the ellipse is assumed known and small, as would be the case if the configuration involved two materials with slightly different electrostatic properties. The problem is to determine the five parameters of the ellipse, given some number of voltage measurements. Two of us (GC and FS) previously studied a version of this problem in which the two electrodes in each pair are co-located (dipole data) \cite{chen-santosa-titi}. In this paper, we extend our previous developments to the more experimentally practical case of separated electrodes.

Even though the model configuration is heavily constrained, some of the parameters (notably the ellipticity and orientation of the anomaly) may be poorly determined. This is particularly the case for arbitrarily chosen electrode positions. To develop a well-posed question for this problem, we use the Jacobian determinant of the model-data relation as a measure of well-posedness, to be optimized as a function of the experiment design parameters (the electrode positions, in this case). This criterion is closely related to so-called \emph{D-optimality}, one of several techniques described in the literature on optimal design of experiments \cite{chaloner-verdinelli,pukelsheim}. In a numerical illustration, we find that optimal placement of electrodes considerably reduces the sensitivity of the inversion to data noise. In effect, the well-posed question in this instance is simply an instance of the inverse problem, with optimally chosen experimental parameters.

In the section following this introduction, we develop a linearized version of the EIT forward map. For the case considered in this paper, of a small elliptical anomaly in a homogeneous background conductivity, several further simplifications are possible. In the next two sections we formulate the restricted EIT problem so defined as a regularized least squares problem. and describe an approach to optimal placement of electrodes. We apply this approach to a numerical example, and observe a dramatic reduction in data sensitivity and model parameter error. A brief discussion of other possible approaches concludes our presentation.

        



\section{EIT for a small elliptical anomaly}

\begin{figure}
\begin{center}
\begin{tikzpicture}[scale=1.3]
\draw[line width=1pt] (0,0) circle (2);
\draw[->] (0,0) -- (2.5,0);
\draw[->] (0,0) -- (0,2.5);
\draw[dotted] (0,0) -- (1.8794,0.6840);
\draw[dotted] (0,0) -- (0.6840,1.894);
\draw[dotted] (0,0) -- (0.3473,-1.9696);
\draw[dotted] (0,0) -- (1.7321,-1.0000);
\draw (0.8,0) arc[start angle=0, end angle=20, radius=0.8];
\draw (0.7,0) arc[start angle=0, end angle=70, radius=0.7];
\draw (0.6,0) arc[start angle=0, end angle=280, radius=0.6];
\draw (0.5,0) arc[start angle=0, end angle=330, radius=0.5];
\draw[black,fill=red] (1.8794,0.6840) circle (.5ex);
\draw[black,fill=blue] (0.6840,1.8794) circle (.5ex);
\draw (2.2,0.7) node[]{$A_+$};
\draw (0.6,2.2) node[]{$A_-$};
\draw[black,fill=green] (0.3473,-1.9696) circle (.5ex);
\draw[black,fill=yellow] (1.7321,-1.0000) circle (.5ex);
\draw (2.1,-1.0) node[]{$B_-$};
\draw (0.3,-2.3) node[]{$B_+$};

\draw (1.1,0.25) node[]{$\phi_+$};
\draw (0.7,0.65) node[]{$\phi_-$};
\draw (-0.1,-0.8) node[]{$\theta_+$};
\draw (0.4,-0.7) node[]{$\theta_-$};
\draw (-0.8,0.7) node[]{$D$};
\draw (2.5,-0.4) node[]{$x_1$};
\draw (-0.4,2.5) node[]{$x_2$};
\filldraw[cyan, opacity = 0.5, rotate=45] (0.2,1.4) ellipse (0.3 and 0.17); 
\end{tikzpicture}
\end{center}
    \caption{In this schematic, the source-sink pair is $A_+$-$A_-$, i.e., current flows from $A_+$ to $A_-$.
    The voltage measurement is taken at $B_+$-$B_-$, i.e., the voltage difference between $B_+$ and $B_-$. The small elliptical anomaly is indicated by $D$.}
    \label{schematic}
\end{figure}
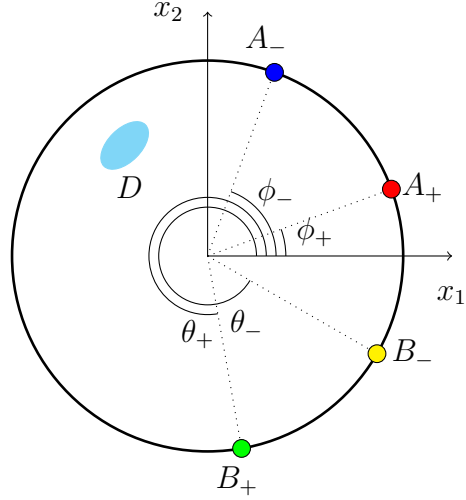
Figure \ref{schematic} depicts a unit disk $\Omega$ in two dimensions with a background conductivity of one. We place a small elliptical anomaly $D$ in $\Omega$. The center of $D$ is $(b_1,b_2)$, the major and minor axes are $(a_1, a_2)$, and the ellipse is oriented at an angle $\xi$ with respect to the $x_1$-axis. The conductivity of the elliptical anomaly is $\sigma=1+\delta$, where $\delta \ll 1$. We abbreviate the parameters of the ellipse as $t = (b_1,b_2,a_1,a_2,\xi)$.

We generate a voltage potential within $\Omega$ by sending unit current through a pair of point electrodes placed on the boundary at $A_{\pm}=(\cos\phi_{\pm},\sin\phi_{\pm})$, see Fig. \ref{schematic}. The voltage potential is approximated as $u=U_0+\delta v$. The background field is given by
\begin{equation}
    U_0(x,y) = \log \frac{(x-\cos\phi_+)^2 + (y-\sin\phi_+)^2}{(x-\cos\phi_-)^2 + (y-\sin\phi_-)^2} .
    \label{background}
\end{equation}
The perturbational field $v$ satisfies
\begin{equation} \label{perturbational}
-\Delta v = \nabla \cdot (\chi_D \; U_0 ), \;\; 
\left. \frac{\partial v}{\partial n}\right|_{x^2+y^2=1} = 0.
\end{equation}
Here $\chi_D$ is the indicator function for domain $D$ and $n$ is the outward normal direction. The relationship between the boundary values of $v$ and $D$ is given by the classical Calderon identity \cite{calderon}
\begin{equation}\label{calderon}
    \int_{\partial B} v \frac{\partial w}{\partial n} ds =
    \int_D \nabla U_0 \cdot \nabla w \; dx dy,
\end{equation}
where $w$ is any harmonic function. Suppose we place electrodes at $B_+: (\cos\theta_+,\sin\theta_-)$ and $B_-: (\cos\theta_-,\sin\theta_-)$ and measure the voltage difference $v(1,\theta_+)-v(1,\theta_-)$ (in polar coordinates).
To obtain a formula relating the perturbational voltage difference above, we choose
\[
w(x,y) = \log \frac{(x-\cos\theta_+)^2 + (y-\sin\theta_+)^2}{(x-\cos\theta_-)^2 + (y-\sin\theta_-)^2} ,
\]
in the Calderon identity \eqref{calderon} and obtain
\begin{equation}
v(1,\theta_+)-v(1,\theta_-) = \int_D \nabla U_0 \cdot \nabla w ;\ dxdy. \label{forward}
\end{equation}

In theory it is entirely possible to use the same locations for both measurement and driving electrode pairs, that is, $A_{\pm}=B_{\pm}$, or equivalently $\theta_{\pm}=\phi_{\pm}$. We will assume this simplification for the remainder of the discussion without delving into the practicality of executing such measurements.

As we assume that the ellipse $D$ is small, a further approximation of \eqref{forward} is possible. Write
\[
P(x,y;\phi_+,\phi_-) =
\nabla U_0(x,y;\phi_+,\phi_-) \cdot \nabla w(x,y;\phi_+,\phi_-),
\]
where we have explicitly denoted the dependence of $U_0$ on electrode locations $(\phi_+,\phi_-)$. Taylor-expand $P$ to second order at the ellipse center $(b_1,b_2)$
\begin{align*}
P(x,y;\phi_+,\phi_-) & =
P(b_1,b_2;\phi_+,\phi_-) \\
& + P_x (b_1,b_2;\phi_+,\phi_-)(x-b_1) +
P_y (b_1,b_2;\phi_+,\phi_-,)(y-b_2)\\ 
& + P_{xx} (b_1,b_2;\phi_+,\phi_-)\frac{(x-b_1)^2}{2} 
+ P_{xy} (b_1,b_2;\phi_+,\phi_-)(x-b_1)(y-b_2) \\ 
&+ P_{xx} (b_1,b_2;\phi_+,\phi_-)\frac{(y-b_2)^2}{2} + \dots.
\end{align*}
Next, attach a coordinate system $(x',y')$ whose origin is at $(b_1,b_2)$, and the axis $x'$ is along the major axis of the ellipse. The change of coordinate is given by
\[
\left[ \begin{matrix}
    x - b_1 \\ y - b_2 
       \end{matrix} \right] =
    \left[ \begin{matrix}
 \cos \xi & -\sin \xi \\
 -\sin \xi & \cos \xi \end{matrix} \right] 
 \left[ \begin{matrix}
     x' \\ y'
 \end{matrix} \right].
\]   
Calling the matrix $R$, we see that the Taylor expansion can be rewritten as
\begin{align}
P(x,y;\phi_+,\phi_-)  =
P(b_1,b_2;\phi_+,\phi_-) 
 &+ [P_x,P_y] R \left[ \begin{matrix}
     x' \\ y'
 \end{matrix} \right] \nonumber  \\
 &+ \frac{1}{2} [x',y'] R^T \left[ \begin{matrix}
    P_{xx} & P_{xy} \\
    P_{xy} & P_{yy}
\end{matrix} \right]
R \left[ \begin{matrix}
     x' \\ y'
 \end{matrix} \right]
 + \dots. \label{taylor}
 \end{align}
 Integrating the expression over $D$, now a standard ellipse in the $x'y'$ coordinates, we get for \eqref{forward}
 \begin{align}
    v(1,\phi_+)-v(1,\phi_-) =
     \pi a_1 a_2 P(b_1,b_2;&\phi_+,\phi_-) + \frac{1}{2} \pi a_1^3 a_2 [\cos\xi,\sin\xi] \left[ \begin{matrix}
    P_{xx} & P_{xy} \\
    P_{xy} & P_{yy}
\end{matrix} \right]
\left[ \begin{matrix}
    \cos\xi \\ \sin\xi
 \end{matrix} \right]
     \label{forward2} \nonumber  \\
     +& \frac{1}{2} \pi a_1 a_2^3 [-\sin\xi,\cos\xi] \left[ \begin{matrix}
    P_{xx} & P_{xy} \\
    P_{xy} & P_{yy}
\end{matrix} \right]
\left[ \begin{matrix}
    -\sin\xi \\ \cos\xi
 \end{matrix} \right]  . 
 \end{align}
 after dropping higher order terms (the next order being $O(a_1^5)$).

\theoremstyle{remark} 
\newtheorem*{remark}{Remark}

Considering the form of the right-hand side of \eqref{forward2}, we introduce new unknowns $A = \pi a_1 a_2$ (ellipse area) and $r=\frac{a_1}{a_2}$ (aspect ratio).
Then the right-hand side of equation \eqref{forward2} can be re-written as
\begin{align}
v(1,\phi_+)-v(1,\phi_-)     = A P(b_1,b_2;\phi_+,\phi_-) + & \frac{1}{2} \frac{A^2 r}{\pi} [\cos\xi,\sin\xi] \left[ \begin{matrix}
    P_{xx} & P_{xy} \\
    P_{xy} & P_{yy}
\end{matrix} \right]
\left[ \begin{matrix}
    \cos\xi \\ \sin\xi
 \end{matrix} \right]
     \label{forward4}  \\
     +& \frac{1}{2} \frac{A^2}{\pi r} [-\sin\xi,\cos\xi] \left[ \begin{matrix}
    P_{xx} & P_{xy} \\
    P_{xy} & P_{yy}
\end{matrix} \right]
\left[ \begin{matrix}
    -\sin\xi \\ \cos\xi
 \end{matrix} \right] .  \nonumber
 \end{align}

We note the `reciprocity' in the linearized forward map: the voltage drop between electrodes at angles $\phi_+$ and $\phi_-$ due to unit current flowing between these electrodes remains the same if $\phi_+$ and $\phi_-$ are swapped. 
This property implies that we can only obtain three independent measurements from three electrodes, and six from four electrodes. Given that we have five unknowns, it seems likely that we must use at least four electrodes to obtain enough measurements to determine all of the ellipse parameters. 


\section{A Regularized Least Squares Formulation}

In view of the preceding observations, we will parametrize the forward map for our version of the EIT problem by a 4-tuple of angles $\Phi=(\phi_1,\phi_2,\phi_3,\phi_4)$. The model parameter vector is $t=(b_1,b_2,A,r,\xi)$. The forward (predicted data) map is then
\[
{\mathcal F}(t,\Phi) = \left[ \begin{matrix}
    v(1,\phi_1)-v(1,\phi_2) \\
    v(1,\phi_1)-v(1,\phi_3) \\
    v(1,\phi_1)-v(1,\phi_4) \\
    v(1,\phi_2)-v(1,\phi_3) \\
    v(1,\phi_2)-v(1,\phi_4) \\
    v(1,\phi_3)-v(1,\phi_4)     
\end{matrix} \right] ,
\]
in which the voltage perturbation $v$ is computed via formula \eqref{forward4}, using the parameters in the vector $t$.
Then the inverse problem consists in choosing the parameter vector $t$ so that the predicted and observed data match, at least approximately:
\begin{equation} \label{IP0}
\mathcal{F}(t,\Phi) \approx g(\Phi).
\end{equation}

Recall that we have assumed that the elliptical anomaly $D$ is small, hence so is its area $A$. Examination of the expression on the right-hand side of equation \eqref{forward4} reveals that the aspect ratio $r$ and orientation angle $\xi$ are scaled by $A^2$ in their effect on the forward map output. Therefore we expect these parameters to be more poorly determined than the others (location coordinates $b_1,b_2$ and area $A$). In formulating a least squares version of the inverse problem \eqref{IP0}, we include a Tihonov regularization term biasing $r, \xi$ towards prior values $r_{\rm prior},\xi_{\rm prior}$. We presume that the data is a value of the forward map applied to a "ground truth" parameter vector $t_0$, corrupted by additive noise $n$:
\[
g(\Phi) = {\mathcal F}(t_0,\Phi) + n,
\]
We also presume that the noise vector $n$ satisfies $\|n\| \leq \epsilon \|g(\Phi)\|$, where $\epsilon$ is known. We pose our version of the EIT inverse problem in least squares form as 
\begin{equation} \label{minimization}
    t_* = \mbox{arg} \; \min_t  \| {\mathcal F}(t,\Phi) - g(\Phi) \|^2 + \lambda \|R(t-t_{\rm prior})\|^2
\end{equation}
where $R = \mbox{diag }(0, 0, 0, 1, 1)$ and $t_{\rm prior}=(0,0,0,r_{\rm prior}, \xi_{\rm prior})$. Tihonov regularization has been studied extensively \cite{hansen1,hansen2,gockenbach}. The additional penalty weight $\lambda$ must be chosen somehow. Several approaches to this choice have been suggested, including some that do not require an assertion of the noise level, such as the L-curve and Generalized Cross-Validation \cite{hansen1,hansen2}. Since we assert a known value (or bound) for $\|n\|$, we use the Morozov discrepancy principle \cite{morozov}, that is, that $\lambda$ should be chosen to make the residual norm at the optimum $\|{\mathcal F}(t_*,\Phi) - g(\Phi)\|$ approximately the same as the asserted data noise level $\epsilon\|g(\Phi)\|$. Since the residual norm increases with $\lambda$, a simple bisection iteration accomplishes this task, at the cost of repeated computations of $t_*$ via equation \eqref{minimization}.

\section{Optimal experiment design}
Up to this point, we have treated the electrode locations ($\Phi$) as passive parameters in the estimation of the model parameters $t$. It is natural to ask whether there are choices of electrode location that minimize the sensitivity of the optimal estimate $t_*$ (equation \eqref{minimization} to data noise. 
To make an optimal etimate of $\Phi$ in this sense, we begin with an \emph{a priori} choice of $\Phi$, denoted $\Phi_*$, and estimate $t_*$ by solving the minimization problem \eqref{minimization} with this choice $\Phi=\Phi_*$ of electrode locations. This process includes determination of a value of the penalty weight $\lambda_*$ via the discrepancy principle and knowledge of the noise level $\epsilon$. 

We presume that the effect of a model perturbation $t_* \rightarrow t_*+p$ is adequately approximated by the linearization of $\mathcal{F}$ at $t_*$, and replace the problem \eqref{minimization} locally by the linear least squares problem
\begin{equation}
\label{linimization}
\min_p \| \mathcal{J}(t_*,\Phi_*) p + \mathcal{F}(t_*,\Phi_*) - g(\Phi_*)\|^2 + \lambda_* \| R p \|^2 ,
\end{equation}
in which $\mathcal{J}$ is the Jacobian of $\mathcal{F}$ (with respect to $t$).
The literature on optimal experimental design offers several criteria \cite{chaloner-verdinelli, pukelsheim} for choosing optimal parameters $\Phi$ for linear least squares problems such as \eqref{linimization}. Two of these criteria are:
\begin{description}[leftmargin=2cm, labelindent=0.5cm]
 \item[\it D-optimal:] Maximize the determinant of $\left[ \, \mathcal{J}(t_*,\Phi)^T \mathcal{J}(t_*,\Phi) + \lambda_* R^T R \, \right]$;
 \item[\it E-optimal:] Maximize the minimum eigenvalue of $\left[ \, \mathcal{J}(t_*,\Phi)^T \mathcal{J}(t_*,\Phi) + \lambda_* R^T R \, \right]$.
\end{description}


\noindent We adopt the D-optimality criterion: that is, we choose $\Phi$ as the solution of
\begin{equation} \label{OED}
\Phi_o = \mbox{arg} \; \max_\Phi \; \det \left[ \, \mathcal{J}(t_*,\Phi)^T \mathcal{J}(t_*,\Phi) + \lambda_* R^T R \, \right] .
\end{equation}
Having computed $\Phi_o$, we compute an optimized estimate $t_o$ by solving the problem \eqref{minimization} again, with $\Phi=\Phi_o$, including an update of the penalty parameter $\lambda_o$.

\section{A numerical example}

We start with a ``ground truth" elliptical anomaly with area $A=0.025$, center at $(b_1=0.452, b_2=-0.165)$, aspect ratio $r=2.32$, and orientation angle $\xi = 0.864$, that is, $t_0 = (0.025,0.452,-0.165,2.323,0.864)$.
In the first experiment, we place the electrodes at $\Phi_*=(0,\pi/2,\pi,3\pi/2)$, and generated data $g(\Phi_*)$ via 
\eqref{forward} and addition of random noise with relative $l_2$ noise level $\epsilon=0.01$. For the regularization term in the least squares problem  \eqref{minimization}, we used prior values $r_{\rm prior}=1$ and $\xi_{\rm prior}=0$. Bisection applied to the discrepancy principle yielded the value $\lambda_* = 0.000594$ for the penalty weight. We obtained the estimated minimizer $t_*=(0.0258,  0.459, -0.153,  0.813, -0.0145) $ via BFGS iteration \cite{NocedalWright}. 

Next, we use the obtained $t_*$ to find the optimal experiment design $\Phi_o$. 
To solve the maximization problem \eqref{OED}, we used Bayesian Optimization \cite{garnett}, as there appeared to exist multiple maxima. The optimized locations $\Phi_o$ are shown in Figure \ref{optimal estimate}. The modified least squares problem \eqref{minimization} with $\Phi=\Phi_o$ produces $\lambda_o=0.001028.$ and $t_o=(0.0251,  0.449, -0.166,  0.444, -0.693)$. Recalling that the aspect ratio $r$ can be swapped with $1/r$ if we also rotate $\xi$ by $\pi/2$, we note that this estimate is very close to the ground truth values. These results are visualized in Figure 3 where one can clearly see the differences in the two recoveries.
\begin{figure}
\centering
\includegraphics[width=3.5in]{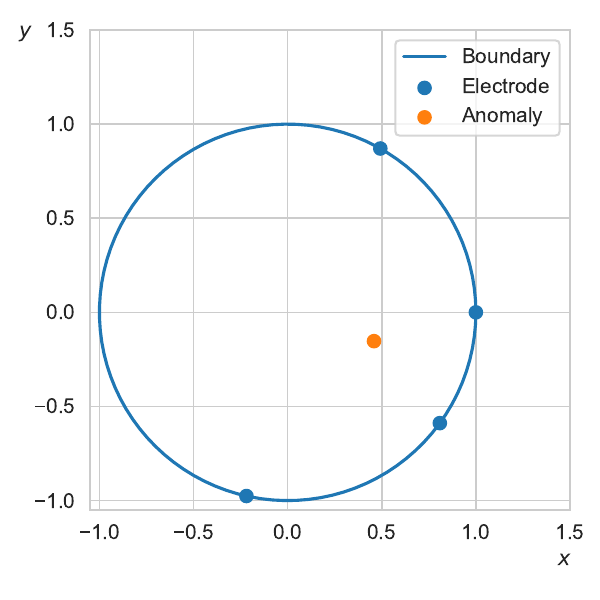}
\caption{The optimal electrodes found by maximizing the determinant of $\left[ \, \mathcal{J}(\Phi)^T \mathcal{J}(\Phi) + \lambda R^T R \, \right]$. The electrodes and the center of estimated anomaly center are indicated. }
    \label{optimal electrodes}
\end{figure}
\begin{figure}
\centering
{\includegraphics[width=2.5in]{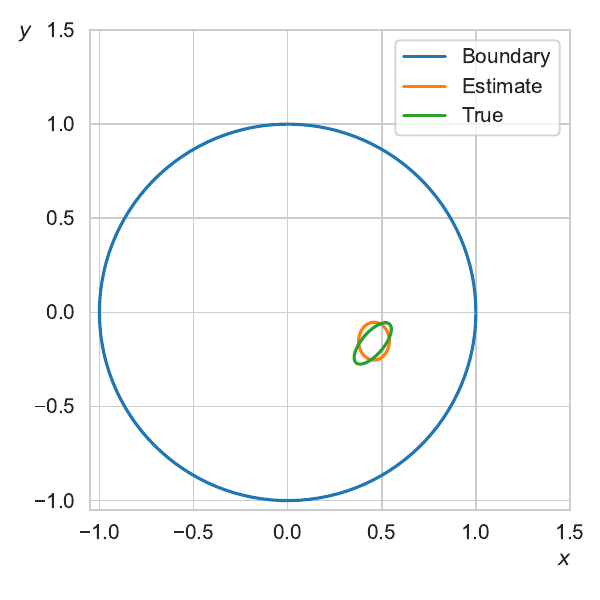}
\includegraphics[width=2.5in]{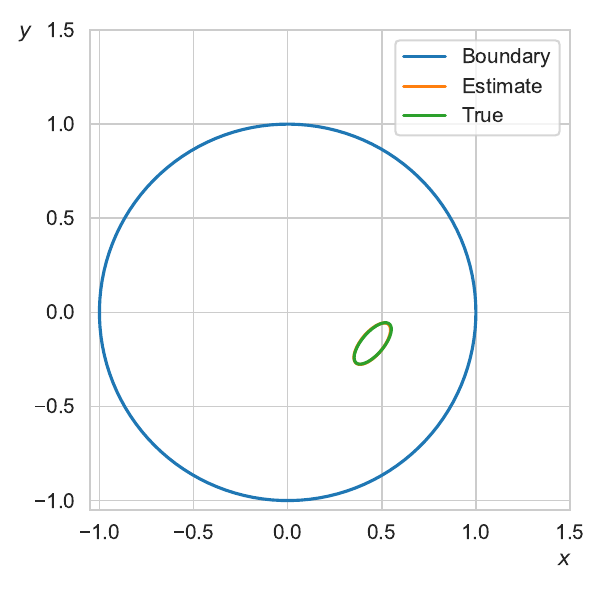}}
\caption{Comparison of the recoveries from two electrode placements with the same relative noise in the data. (Left) The recovery from electrodes placed at $(0,\pi/2,\pi,3\pi/2)$ and (Right) the recovery from optimally placed electrodes (see Figure \ref{optimal electrodes}). Note that in the case of optimally placed electrodes regularization is below the threshold of sensitivity and is basically ignored. }
\label{optimal estimate}
\end{figure}

\begin{figure}
\centering
\includegraphics[width=6.5in]{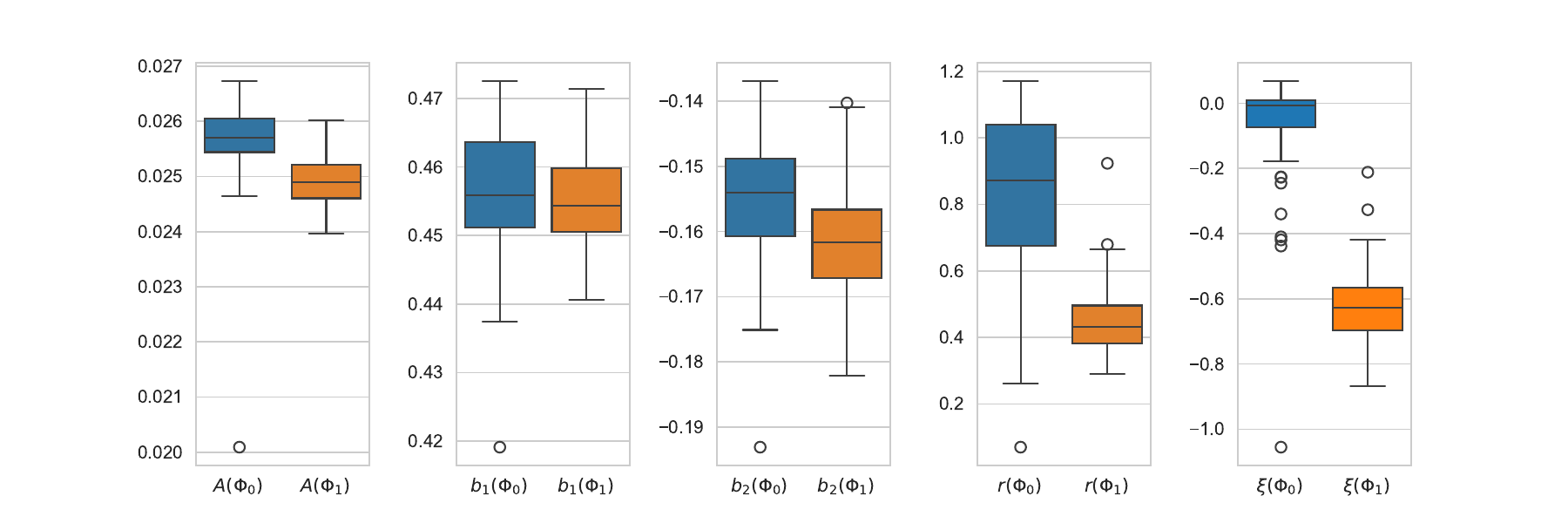}
\caption{These box plots summarizes the 100 inversions using electrodes at initial locations  $(0,\pi/2,\pi,3\pi/2)$ (in blue) and those using optimal electrode locations (in orange - see Figure \ref{optimal electrodes}). Notice how much closer the means of the ellipse parameters for the estimates with optimal electrode placement are to the ground truth values ($A=0.025, b_1=0.452, b_2=-0.165, r=1/2.323, \xi=0.864-\pi/2$). The result suggests that by optimizing the electrode locations, the problem becomes much more well-posed.}
\label{boxes}
\end{figure}

To explore the accuracy and stability induced by optimal electrode placement, we solved \eqref{minimization} 100 times (each time with a new realization of noise at a fix variance and an upper bound) for the two electrode configurations, using the already-computed values $\lambda_*$ and $\lambda_o$ of the penalty weight. The results are summarized in the box plots displayed in Figure \ref{boxes}. The sample means for the first three parameters are reasonably accurate approximations to the``ground truth" in both cases, though more accurate and with slighly smaller variance for the optimized electrode placements. For the last two parameters, the story is completely different: these parameters ($r$ and $\xi$) are close to their prior values for the first (evenly spaced) electrode configuration. That is, in this case, estimation of these parameters is an ill-posed problem, and regularization dominates the estimate. For the optimal electrode placement, these two parameter means are quite close to ``ground truth'', suggesting that optimal design has dramatically increased the sensitivity of the data prediction to these parameters, effectively rendering the inversion well-posed. 

\section{Discussion}

Our approach to optimal experimental design is deterministic: we maximize a function related to the forward map to choose a design. Much of the literature on this topic, however, has a statistical bent. We briefly describe one of these alternative approaches.

An alternate way to formulate well-posed question is to put the inverse problem in the framework of Bayesian inference. This idea is well described in the monograph of Kaipio and Somersalo \cite{kaipio-somersalo}. In this framework, one would start with a prior distribution for the unknown $t$, which we denote here by $P(t)$. The prior distribution embodies the amount of information we have about $t$ prior to collecting data associated with the inverse problem. For example, we can model the distibution as a product of four distributions - one for $A$ (which must be positive), one for $(b_1,b_2)$ which must satisfy $(b_1^2+b_2^2) \leq \rho < 1$, one for $r$ with $r_- \leq r \leq r_+$ and one for $\xi$ with $0\leq \xi <\pi$. All four distributions could be truncated normal distributions with prescribed means and standard deviations.

The key to the Bayesian approach is to obtain a posterior distribution after data is observed, i.e., $P(t\,|\,g; \Phi)$ where $g(\Phi)$ is the measured data. We denote the dependence of the distribution on the experiment $\Phi$. 
Equipped with such the posterior distribution, one can proceed to obtain confidence bounds on the unknown $t$ using Monte Carlo sampling. Due to the lack of sensitivity of the forward map \eqref{forward4} in $r$ and $\xi$, we expect the confidence bounds to be wide in these two parameters, while relatively tight in the other three parameters. One major advantage of the Bayesian approach is that it can be applied to nonlinear problems.

Chaloner and Verdinelli \cite{chaloner-verdinelli} provided a review of Optimal Experiment Design in Bayesian inference. The formulation works for nonlinear problems as shown in Huan and Marzouk \cite{huan-marzouk-2013}. Every experiment $\Phi$ leads to some information gain on the unknown $t$. This information gain can be measured by the Kullback-Leibler divergence.
The KL divergence depends on $g$ and $\Phi$. The expected information gain is the utility of an experiment $\Phi$ and is maximized over possible experiments.
This has been implemented in \cite{chen-santosa-titi} for finding optimal dipole locations to best determine the location and the area of an elliptical anomaly in linearized EIT.

\section{Conclusion}
In the spirit of Sabatier, we have shown that the ill-posed problem of linearized EIT for a small elliptical anomaly can be made well-posed by selective regularization and optimal experiment design. By placing electrodes in such a way as to maximize linear stability, we obtain parameter estimates dominated by the data, rather than partly by the prior (as is the case for non-optimal electrode placement). That is, for this ill-posed problem, the well-posed question arises as a carefully chosen instance of the problem itself.

\end{document}